\documentstyle{article}

\newcommand{\1}{{{\mathchoice {\rm 1\mskip-4mu l} {\rm 1\mskip-4mu l}
{\rm 1\mskip-4.5mu l} {\rm 1\mskip-5mu l}}}}
\newcommand{\R}{{\bf R}}

\newcommand{\Sp}{{\rm Sp}}
\newcommand{\SP}{{\rm Sp}}
\newcommand{\Sl}{{\rm SL}}
\newcommand{\spl}{{\bf sp}}
\newcommand{\SL}{{\rm SL}}
\newcommand{\GL}{{\rm GL}}
\newcommand{\tr}{{\rm tr}}
\newcommand{\diag}{{\rm diag}}

\newcommand{\Rr}{{\cal R}}
\newcommand{\Oo}{{\cal O}}

\newcommand{\Co}{{{\cal C}\!{\it onj}}}

\newcommand{\Hh}{{\bf H}}

\newcommand{\C}{{\bf C}}

\newcommand{\G}{{\rm G}}

\newcommand{\al}{{\alpha}}

\newcommand{\be}{{\beta}}

\newcommand{\om}{{\omega}}
\newcommand{\eps}{{\varepsilon}}
\newcommand{\de}{{\delta}}

\newcommand{\ga}{{\gamma}}

\newcommand{\ka}{{\kappa}}
\newcommand{\la}{{\lambda}}
\newcommand{\si}{{\sigma}}

\newcommand{\Bb}{{\cal B}}
\newcommand{\Cc}{{\cal C}}
\newcommand{\Dd}{{\cal D}}

\newcommand{\Nn}{{\cal N}}

\newcommand{\Pp}{{\cal P}}
\newcommand{\Qq}{{\cal Q}}
\newcommand{\Uu}{{\cal U}}

\newcommand{\Ss}{{\cal S}}

\newcommand{\Ham}{{\rm Ham }}

\newcommand{\Conj}{{\rm Conj}}

\newcommand{\MS}{{\medskip}}

\newcommand{\NI}{{\noindent}}
\newcommand{\proof}[1]{\noindent{\bf Proof#1:\  }}

\newcommand{\QED}{\hfill$\Box$\medskip}

\newtheorem{theorem}{Theorem}[section]

\newtheorem{definition}[theorem]{Definition}

\newtheorem{remark}[theorem]{Remark}
\newtheorem{lemma}[theorem]{Lemma}

\newtheorem{prop}[theorem]{Proposition}
\newtheorem{proposition}[theorem]{Proposition}
\newcommand{\at}{{@}}

\title{Positive paths in the linear symplectic group}
\author{Fran\c{c}ois Lalonde\thanks{Partially supported by NSERC grant
OGP 0092913 
and FCAR grant ER-1199.} \\ Universit\'e du Qu\'ebec \`a Montr\'eal
\\ (flalonde\at math.uqam.ca) \and Dusa McDuff\thanks{Partially supported by
NSF grant DMS 9401443.  Both authors supported by NSERC grant CPG 0163730.} 
\\ State University of New York at Stony Brook \\
(dusa\at math.sunysb.edu)}

\date{August 22, 1995}

\begin{document}

\maketitle

\section{Introduction}

A positive path in the linear symplectic group $\Sp(2n)$ is a  
smooth path
which is everywhere tangent to the positive cone.  
These paths  are generated by negative
definite (time-dependent) quadratic Hamiltonian 
functions on Euclidean space.  A
special case are autonomous positive  paths, which are generated 
by time-independent 
Hamiltonians, and which all lie in the set $\Uu$ of diagonalizable 
matrices with
eigenvalues on the unit circle.  However, as was shown by 
Krein, the
eigenvalues of a general positive path can move off the unit circle.  
In this paper, we extend Krein's theory: we investigate the general
behavior of positive paths which do not encounter the
eigenvalue $1$, showing, for example, that any such path can be extended
to have endpoint with all eigenvalues on the circle. We also
show that in the case $2n=4$ there is a close relation between the index of a
positive path  and the regions of the symplectic group that such a path can
cross. 
Our motivation for studying these paths came from a geometric
squeezing problem~\cite{LM2} in symplectic topology.   However,
they are also of interest in relation to the stability of periodic 
Hamiltonian
systems~\cite{GL} and in the theory of geodesics in Riemannian 
geometry~\cite{B2}.

\subsection*{Main results}

We consider $\R^{2n}$ equipped with the standard (linear) symplectic form
$$
\om(X,Y) = Y^T JX,
$$
where  $J$ is multiplication by $i$ in $\R^{2n}\cong\C^n$.  Thus if 
$e_1,\dots,e_{2n}$ form the standard basis $Je_{2i-1} = e_{2i}$.
The  Lie algebra $\spl(2n)$ of $\G = \SP(2n, \R)$ is the set of all matrices
which satisfy the equation
$$
A^T J + JA = 0.
$$
Hence $A\in \spl(2n)$ if and only
if $JA $ is symmetric. The
symplectic gradient (or Hamiltonian vector
field) $X_H$  of a function $H:\R^{2n}\to \R$ is defined by the identity
$$
\om (X_H,Y) = Y^T J X_H =  dH(Y)\quad \mbox{for all\,}\,Y\in \R^{2n}.
$$
 Given a symmetric matrix $P$, consider the associated function
$$
\Qq (x) = - \frac 12 x^T P x
$$
on $\R^{2n}$.  Then the
calculation $$
 d\Qq_x(Y) =  - Y^T P x = Y^T J(J Px)
$$
shows that the symplectic gradient of $\Qq$  is the vector field
$$
X_\Qq(x)  = JP(x),\quad\mbox{at}\quad x\in \R^{2n},
$$
which integrates to the linear flow
$
x \mapsto  e^{JPt}x.
$
We call these paths $\{e^{JPt}\}_{t\in [0,1]}$ {\bf  autonomous}, since they 
are generated by autonomous (i.e. time-independent) Hamiltonians.

We are interested in paths $A_t\in\G$ which are {\bf positive} in the
sense that their tangent vector at time $t$ has the form $JP_tA_t$ for
some symmetric, positive definite matrix $P_t$.  Thus they are
generated by a family of time-dependent quadratic
Hamiltonians\footnote { It is somewhat unfortunate that our sign
conventions imply that $\Qq_t$ is actually negative definite.  }  of
the form $\Qq_t = -\frac 12 x^T P_t x$.  These paths are everywhere
tangent to the {\bf positive cone field} $\Pp$ on $\G$ given by $$
\cup_{A\in \G}\Pp_A = \cup_{A\in \G} \{JPA: P=P^T, P>0\} \;\subset
\;\cup_A T_A\G\; = \; TG.  $$ In particular, an autonomous path
$e^{JPt}$ is positive if $P$ is.

Note that an autonomous path
$e^{JPt}$ always belongs to the set $\Uu$ of  elements in $\G$ which are 
diagonalizable
(over $\C$) and have spectrum on the unit circle.  The following result shows
that this is not  true for all positive paths.

\begin{prop}\label{prop:pi0}  Any two elements in $\G$ may be joined by a
positive path. \end{prop}

This is a special case of a general result in control theory known as Lobry's
theorem: see Lobry~\cite{LOB}, Sussmann~\cite{SU} and 
Grasse--Sussmann~\cite{GS}.  We give a  proof  here in
\S\ref{ss:ge3} as a byproduct of our other results.

The structure of positive paths $A_{t\in [0,1]}$ in $\G$ was studied
by Krein in a series of papers during the first half of the
50's. Later Gelfand--Lidskii~\cite{GL} took up the study of the
topological properties of general paths in connection with the
stability theory of periodic linear flows.  Almost simultaneously, and
surely independently, Bott~\cite{B1,B2,B3} studied positive flows in
the complexified linear group in connection with the geometry of
closed geodesics. As a result of these works, the structure of
positive paths is well-understood provided that the spectrum of $A_t$
remains on the unit circle. However, we need to investigate their
behavior outside this circle. Surprisingly, a thorough study of the
behavior of these paths outside the unit circle does not seem to
exist.  The aim of this article is to undertake such a study, at least
in the case of generic positive paths, that is to say those paths
which meet only the codimension $0$ and codimension $1$ strata of the
matrix group.  For the convenience of the reader we will develop the
theory from scratch, even though some of our first results are
well-known.

Our main result is motivated by the geometric application
in~\cite{LM2}.  Before stating this result, we introduce some
notation.  As above, $\Uu$ denotes the set of elements in $\G$ which
are diagonalizable and have spectrum on the unit circle (unless
indicated to the contrary, ``diagonalizable'' is considered over
$\C$).  Further, $$
\Ss_1 = \{A\in \G: \det(A - \1) = 0\}
$$ is the set of elements with at least one eigenvalue equal to $
1$. A path $\ga = A_{t\in[0,1]}$ in $\G$ which starts at the identity
$A_0 = \1$ and is such that $A_t\in \G - \Ss_1$ for all $t > 0$ will
be called {\bf short}. Equivalently, it has Conley-Zehnder index equal
to $0$.\footnote { The Conley--Zehnder index (or Maslov index) of a
path measures how many times it goes through the eigenvalue $1$ and
characterizes its homotopy class: see Gelfand--Lidskii~\cite{GL},
Ekeland~\cite{EK1,EK}, Robbin--Salamon~/cite{ROSA1}. }  Moreover, an
element $A\in
\G$ is called {\bf generic} if all its eigenvalues (real or complex) have
multiplicity $1$. In the next theorem, all paths begin at $\1$.

\begin{theorem}\label{thm:pospath}
\begin{description}
\item[(i)]  An 
element of $\G - \Ss_1$ is the endpoint of  a short positive 
path if and only if it has an even
number of real eigenvalues $\la$ with $\la > 1$.  
\item[(ii)] Any short
positive path  may be extended to a short positive path  with
endpoint in $\Uu$.
\item[(iii)]  The space of short positive paths with endpoint in $\Uu$
is path-connected. \end{description}
\end{theorem}

We can interpret these results in terms of a \lq\lq conical"
subRiemannian geometry.  SubRiemannian geometry is usually a study of
paths in a manifold $M$ which are everywhere tangent to some
distribution $\xi$ in the tangent bundle $TM$.  For example, there is
recent interest in the case when $\xi$ is a contact structure:
see~\cite{Gr}.  One could generalize this to the consideration of
paths which are everywhere tangent to some fixed convex conical
neighborhood $\Nn_\eps(\xi)$ of $\xi$.  The geometry of positive paths
more or less fits into this context: one just has to replace the
distribution $\xi$ by a distribution of rays.  Thus we can consider
the positive cone to be a neighborhood of the ray generated by the
right invariant vector field $JA\in T_A\G$ corresponding to the
element $J$ in the Lie algebra.  Proposition~\ref{prop:pi0} shows that
$\G$ remains path-connected in this geometry.  Next, one might try to
understand what happens to the fundamental group.  More precisely, let
us define the {\bf positive fundamental group} $
\pi_{1,{\rm pos}}\G 
$ to be the semigroup generated by positive loops based at $\1$, where
 two such loops are considered equivalent if they may be joined by a
 smooth family of positive loops.  It would be interesting to
 calculate this semigroup.  In particular it is not at present clear
 whether or not the obvious map $$
\pi_{1,{\rm pos}}\G \;\to\; \pi_1(\G)
$$
is injective in general.   (When $n = 1$ this is an easy consequence of 
Proposition~\ref{prop:lift} and Lemma~\ref{le:pos1}.)

The proof of Theorem~\ref{thm:pospath} is based on studying the
relation between the positive cone and the fibers of the projection
$\pi: \G \to \Co$, where $\Co$ is the space of conjugacy classes of
elements in $\G$ with the quotient (non-Hausdorff) topology.  The
space $\Co$ has a natural stratification coming from the
stratification of $\G$.  Recall that a symplectic matrix $A$ is
similar to its transpose inverse $(A^T)^{-1}$ (in fact, $(A^T)^{-1} =
-JAJ$).  Thus its eigenvalues occur either in pairs $\la,
\bar\la \in S^1$ and $\la, 1/\la \in \R - 0$, or in complex 
quadruplets $\la, \bar\la, 1/\la,
1/\bar\la$. In particular, the eigenvalues $\pm 1$ always occur with 
even multiplicity.  Roughly speaking,  the open strata of $\G$ consist of
generic  elements (ie
diagonalizable elements with all eigenvalues of multiplicity $1$), and 
 the
codimension $1$ strata consist of non-diagonalizable elements 
which have one pair
of eigenvalues of multiplicity $2$ lying either on $\R - \{\pm 1,0\}$ 
or on $S^1 - \{\pm
1\}$, or which have one pair of eigenvalues equal to $\pm 1$: 
see~\cite{AG} and \S3 below.   Note that a conjugacy class which 
lies in an open
stratum is a submanifold of ${\rm GL}(2n, \R)$ of codimension $n$ 
since each pair 
of distinct eigenvalues (in
$\R \cup S^1$) has one degree of freedom, 
and each quadruplet
has two.
 
    We write $\G_0$ for the set of all generic elements in $\G$, 
$\G_1$ for the union of all codimension $1$ strata, and similarly  $\Co_0,
\Co_1$ for their projections in $\Co$. We shall
consider in detail only generic paths in $\Co$.  By definition, these
intersect all lower strata transversally.  Hence the elements of these
paths lie in $\Co_0$ for all except the finite number of times at
which they cross $\Co_1$. Note that the codimension of a stratum in
$\Co$ always refer to the codimension of its lift in $\G$.

The main ingredient of the proof of the theorem is to characterize the
(generic) paths in $\Co$ which lift to positive paths in $\G$.  We
shall see in Proposition~\ref{prop:pos2} that the only significant
restriction on the path comes from the way that eigenvalues enter,
move around and leave $S^1$.  We will not state the general result
here since it needs a certain amount of notation.  However some of the
essential features are present in the case $n = 1$ which is much
easier to describe.  In this case $\Co$ is the union of the circle
with the intervals $(-\infty, -1]\cup [1,\infty)$ with the usual
topology, except that the points $\pm 1$ are each tripled.  To see
this, first consider matrices with an eigenvalue pair $\la, \bar\la
\in S^1$.  It is possible to distinguish between these eigenvalues and
hence to label the pair with an element in $S^1$.  In higher
dimensions this is accomplished by the {\it splitting number}
described in \S2.  Here it suffices to note that the rotations through
angles $\theta$ and $2\pi - \theta$ are not conjugate in $\SL(2,\R)$,
the only conjugating matrices being reflections.  We also claim that
the point $1$ occurs with three flavors: plain $\1$ (the conjugacy
class of the identity) which has codimension $3$, and two nilpotent
classes $\1^\pm \in \Co_1$.  A similar statement holds at $-1$.  As
explained in Lemma~\ref{le:pos1}, a positive path has to project to a
path in $\Co$ which goes round the circle anticlockwise, but it can
move along the real axes in either direction.  Moreover, the
projection of a generic positive path has to leave the circle at
$(-\1)^-$ or $\1^-$ and then, after wandering around $(-\infty, -1)
\cup (1,\infty)$ for a while, enter the circle again at $(-\1)^+$ or
$\1^+$ and continue on its way.

   The next step is to carry out a similar analysis in the case $n=2$.
We shall see that the restrictions in the behavior of positive paths
in this case all stem from Krein's lemma which says that simple
eigenvalues on $S^1$ with positive splitting number must flow
anti-clockwise round the circle.  This gives rise to restrictions in
the way in which four eigenvalues on the unit circle can move to a
quadruplet outside the unit circle (the Krein--Bott bifurcation).  As
a consequence, we show in Proposition~\ref{prop:restr} that there are
restrictions on the regions in $\Co$ which can be visited by a generic
positive path of bounded Conley-Zehnder index.

   In the last section of the paper, we will present a brief survey of 
the use of positive paths in 
stability theory, in Riemannian geometry and in Hofer geometry, and 
will give some new applications.

We wish to thank Nancy Hingston for describing to us her ideas about 
perambulations in
the symplectic group, which were a direct inspiration for our work, and Hector
Sussmann for explaining to us some basic results in  control theory.

\section{Basic facts}
\subsection{General results on positive paths}

First, two elementary lemmas.  A piecewise smooth path is said to be 
piecewise positive if
the tangent vectors along each smooth segment (including those at 
the endpoints) are in the positive cone. 

\begin{lemma}\label{le:open} \begin{description}
\item[(i)] The set of
positive paths is open in the $C^1$-topology.  
\item[(ii)] Any piecewise positive
path may be $C^0$-approximated by a positive path. 
\item[(iii)]  Let $\{A_t\}, \{B_t\}$ be two positive paths with
the same initial point $A_0 = B_0$.  Then, given $\eps > 0$, there is 
$\de > 0$
and a $C^0$-small deformation of $\{A_t\}$ to a positive path $\{A_t'\}$ which
coincides with $\{B_t\}$ for $t \in [0,\de]$ and with $\{A_t\}$ for $t >
\eps$. 
\end{description}\end{lemma}

\proof{} The first statement follows from the openness of the positive
cone,  the second from its convexity and the third from a
combination of these two facts.  \QED

Note that statement (ii) above is a special case of a much more general
result in control theory: see, for example, Sussmann~\cite{SU} and
Grasse--Sussmann~\cite{GS}.

 \begin{lemma}\label{le:gen} 
\begin{description} \item[(i)]  Any positive path  can be
extended to a positive  path which
ends at a generic element.
\item[(ii)]  Any positive path $\{A_t\}$ starting at $A_0=\1$ 
can be perturbed fixing
$A_0$ so that $\{A_t\}$ is generic.  Moreover, if $A_1$ is 
generic we may fix that
also during the perturbation. \end{description}\end{lemma}
\proof{}  This follows immediately from openness.\QED

We continue with some simple remarks about the relation between positive
paths in $\G$ and their projections (also called positive paths) in $\Co$.

\begin{lemma}\label{le:conj0}  
\begin{description}
\item[(i)] The positive cone is invariant under conjugation.
\item[(ii)]  The 
conjugate $\{X^{-1}A_t X\}$ of  any positive path $\{A_t\}$ is positive.
\item[(iii)] There is a positive path from $\1$ to $A$ in $ \G$ if and
only if there is a positive path in $\Co$ from $\1$ to $\pi(A)$.
\end{description}
\end{lemma}
\proof{}  (i)  The positive vector $JPA$ at $A$  is taken by conjugation to the
positive
vector $X^{-1}JPAX  = X^{-1}JPX X^{-1}AX = J(X^TPX)B$ at $B = X^{-1}AX$.

\NI
(ii) follows immediately from (i).

\NI
(iii) Since the \lq\lq only if" statement is obvious, we consider the
\lq\lq if" statement.  By definition, a positive path from $\pi(\1)$
to $\pi(A)$ lifts to a positive path from $\1$ to $Y^{-1}AY$ for some
$ Y\in G$.  Now conjugate this path by $Y^{-1}$ and use (ii).\QED

\NI
{\bf Warning} Given a positive path $\{A_t\}$ it is {\it not} true
that all paths of the form $\{X_t^{-1}A_t X_t\}$ are positive.  One
can easily construct counterexamples in $\SL(2,\R)$ using the methods
of \S 3.

\begin{prop}\label{prop:lift}   
Let $\{A_t\}$ in $\G$ be a generic positive path joining two generic
points $A_0, A_1$.  Then the set of positive paths in $\G$ which lift
$\ga = \{\pi(A_t)\}_{t\in [0,1]}\subset\Co$ is path-connected.
Moreover, the set of these paths with one fixed endpoint is also
path-connected. Finally, if $A_0=\1$, the set of these paths with both
endpoints fixed is path-connected.
\end{prop} 

\proof{}  
Let $\{B_t\}$ be another path which lifts $\ga = \{\pi(A_t)\}$.  By
genericity we may suppose that $\{A_t\}$ and $\{B_t\}$ are disjoint
embedded paths. Moreover, when $n
\ge 2$, the projection in $\Co$ of the 
positive cone at any point $X \in \G_0 \cup
\G_1$ of a  generic path is always at least $2$-dimensional (see
Proposition~\ref{prop:pos2}).  Hence we may assume that the paths
$\{A_t\}$ and $\{B_t\}$ are never tangent to the fiber of the
projection $\G \to \Co$.  Our first aim is to define positive
piecewise smooth vector fields $\xi_A, \xi_B$ that are tangent to
$\pi^{-1}(\ga)$ and which extend the tangent vector fields to the
paths $\{A_t\},
\{B_t\}$.  
To make this easier we first normalise the path $B_t$ near the finite
number of times $t_1 <
\dots < t_k$  at which  $A_t$ crosses a 
codimension $1$ stratum.  If $T$ is one of these times $t_i$, there is
a matrix $X_T$ such that $B_T = X^{-1}_T A_T X_T$, and by
Lemma~\ref{le:open} (iii) we may suppose that $B_t = X^{-1}_T A_t X_T$
for $t$ near $T$.

On the other hand, since the set of tangent vectors $Y \in T\G$ which
projects to the tangent vector $\pi_{\ast}(\dot{A}_T)$ is an affine
subbundle of the restriction $T\G \mid_{\pi^{-1}(\pi(A_T))}$, we can
find a section of this subbundle whose value at $A_T$ coincides with
$\dot{A}_T$. This defines $\xi_A$ over neighbourhoods of all times
$t_i$, which we choose small enough so that $\xi_A$ is positive. We
then define $\xi_B$ over each such neighbourhood by the adjoint map
$X^{-1}_T \xi_A X_T$.

At all other times $A_t$ is generic and the map $t\mapsto \pi(A_t)$ is
an immersion. Therefore, the set $$ Y_i = \{\pi^{-1}(\pi(A_t)): t_i\le
t\le t_{i+1}\} $$ 
is a submanifold.  We now extend $\xi_A, \xi_B$ to
the whole of $Y_i$ by choosing sections of the intersection of the
positive cone field with $TY_i$ which project to the tangent vector
field of $\ga$.  This is possible because this intersection is an open
convex cone.  Observe that $\pi$ projects all integral curves of
$\xi_A$ and $\xi_B$ to $\ga$.  Indeed all integral curves of the
vector field $s\xi_A + (1-s)\xi_B$ (for $s\in[01]$) are positive paths
which project to $\ga$.  Clearly there is a path of such curves
joining $\{A_t\}$ to $\{B_t\}$.

It remains to check the statement about the endpoints.  We will
assume, again using Lemma~\ref{le:open}, that if $\{A_t\}$ and
$\{B_t\}$ have the same endpoint they agree for $t$ near that
endpoint.  Then, if $B_0 = A_0$ the above construction gives a family
of paths starting at $A_0$.  Therefore the result for fixed initial
endpoint is obvious.  To get fixed final endpoint one can use a
similar argument applied to the reversed (negative) paths
$\{A_{1-t}\},\{ B_{1-t}\}$.  Finally consider the case of two fixed
endpoints, with $A_0=B_0=\1$.  The endpoints $A_1^\la$ of the family
of paths constructed above are all conjugate to $A_1$ and equal $A_1$
at $\la = 0,1$.  Because $A_1$ is generic, there is a smooth family of
elements $X_\la\in \G$ with $X_0, X_1 = \1$ such that $X_\la^{-1}
A_1^\la X_\la = A_1$.  The desired family of paths is then
$\{X_\la^{-1} A_t^\la X_\la\}_{t\in [0,1]}$: these are positive by
Lemma~\ref{le:conj0}. \QED

\subsection{Splitting numbers}

This section summarizes what is known as Krein theory, that is the
theory of positive paths in $\Uu$.  The theory is nicely described by
Ekeland in~\cite{EK1,EK} where it is used to develop an index theory
for closed orbits.  Interestingly enough, this analysis is very
closely related to work by Ustilovsky~\cite{UST} on conjugate points
on geodesics in Hofer geometry on the group of Hamiltonian
symplectomorphisms.

It will be convenient to work over $\C$ rather than
$\R$.  We extend $J$ and $\om$ to $\C^{2n} = \R^{2n}\otimes \C$ by
complex linearity, so that $\om(v,w) = w^T Jv$ is complex
bilinear. Let $\langle v, w\rangle$ denote the standard Hermitian inner 
product on $\C
^n$, namely
$$
\langle v, w\rangle = \bar w^T v = \om(Jv,\bar w).
$$
 We will also use  the form $\be$ given in our notations by
$$
\be(v,w) =   -i \om(v,\bar w) = -i \langle Jv, w\rangle = -i \bar w^T Jv.
$$

\begin{lemma}\label{le:be}
\begin{description}
\item[(i)] $\be$ is a nondegenerate Hermitian symmetric form.

\item[(ii)]   $i\be(v,Jv) = \langle v, v\rangle$.

\item[(iii)] An element $A\in \GL(2n,\C)$ belongs to $\G$ if and only if it
preserves $\be$ and is real, i.e.
$$
\overline {Av} = A\bar v,\quad v\in \C^{2n}. 
$$

\item[(iv)] The invariant subspaces $E_\la, E_\mu$ corresponding to eigenvalues
$\la,\mu$ are $\be$-orthogonal unless $\la\bar\mu = 1$.

\item[(v)]  If $v$ is an eigenvector with eigenvalue  $\la\in S^1$ of
multiplicity $1$ then $\be(v,v)\in \R - \{0\}$.
\end{description}\end{lemma}
\proof{}  $\be$ is Hermitian because
$$
 \overline{\be( w, v)} =\overline{-i\om( w, \bar v)} = i{\om(\bar w,v)} =
-i\om(v,\bar w) = \be(v,w).
$$
Statements (ii) and (iii) are obvious.  
To prove (iv), consider $A\in \G = \Sp(2n,\R)$ with eigenvalues
$\la,\mu$ and corresponding eigenvectors $v,w$ and observe that
$$
\be(v,w) = \be(Av,Aw) = \be (\la v, \mu w)  = \la\bar \mu\be(v,w).
$$
The proof when $v,w$ belong to the invariant subspace but are not eigenvectors
is similar (see Ekeland \cite{EK}).
To prove (v), suppose $\la\in S^1$ is an eigenvalue of
multiplicity $1$ with eigenvector $v$.  Then $\be(v,v)\ne 0$, since
$\be$ is non-degenerate and $\be(v,w) = 0$ whenever $w$ belongs to
any other eigenspace of $A$.  
Further, $\be(v,v) \in \R$ by (i). 
\QED

It follows from (v) that simple eigenvalues $\la$ on
$S^1$ may be labelled with a number  $\si(\la)=\pm 1$ called the {\bf
splitting number} chosen so that
$$
\be(v,v) \in  \si(\la) \R^+.
$$
For example, when $n = 1$, the rotation matrix
$$
\rho(t) = \left(\begin{array}{cc}  \cos t & -\sin t\\ \sin t & \cos t
\end{array}\right)
$$
has eigenvectors
$$
v^+ = \left(\begin{array}{c}  1\\-i
\end{array}\right),\quad 
v^- = \left(\begin{array}{c}  1\\i
\end{array}\right)
$$
corresponding to the eigenvectors $e^{it}, e^{-it}$, and 
it is easy to see that $e^{it}$ has splitting number $+1$, while
$e^{-it}$ has splitting number $-1$.   

More generally, if $\la\in S^1$ has
multiplicity $> 1$ we define $\si(\la)$  to be the signature of the form
$\be$ on the corresponding invariant subspace.  Note that $\si(\la) =
-\si(\bar\la)$.  Hence the splitting number of $\pm 1$ is always $0$.  
It is not
hard to check  that the conjugacy class  of an element in $\Uu$  is completely
described by its spectrum together with the corresponding set of splitting
numbers.  The following result is central to our argument.  

\begin{lemma}[Krein]\label{le:bott}  
Under a positive flow simple eigenvalues on
$S^1$ labelled with $+1$ must move anti-clockwise and those labelled
with $-1$ must move clockwise.
\end{lemma}   
\proof{}  We repeat the proof from Ekeland's
book~\cite{EK} for the convenience of the reader.  For all $t$ in some
neighborhood of $t_0$, let $ e^{i\theta_t}$ be a simple eigenvalue of
$A_{t}$ with corresponding eigenvector $x_t\in \C$.  For simplicity
when $t = t_0$ we use the subscript $0$ instead of $t_0$ (writing
$A_0, x_0$ etc), and will denote the derivative of $A_t$ at $t_0$ by
$JPA_0$, where $P$ is positive definite.  It is easy to check that for
any $x$ $$
\langle A_t x, Jx_t\rangle = \langle x, A_t^T Jx_t\rangle
= \langle x, JA_t^{-1}x_t\rangle = \langle  e^{i\theta_t} x, Jx_t\rangle.
$$
Applying this with $x = x_{0}$ and differentiating at $t = t_0$, we find
that 
$$
\langle JP A_0 x_{0}, Jx_{0}\rangle = \left. \langle i \frac{ d\theta_t}{dt} 
e^{i\theta_0}x_0, Jx_0\rangle \right|_{t = t_0},
$$
from which it readily follows that
$$
\left. \frac{ d\theta_t}{dt}  \right|_{t = t_0}= 
\frac{\langle Px_{0}, x_0\rangle}{\be(x_0,x_0)}.
$$
Since the right hand side has the same sign as ${\be(x_0, x_0)}$, the result
follows. \QED

 Observe, however, that not all flows whose eigenvalues move in this
way are positive.  Further, this result shows that there may not be a
short positive path between an arbitrary pair of elements in a given
conjugacy class, even if one allows the path to leave the conjugacy
class.

The next lemma shows that if a (simple) eigenvalue leaves
the circle it must do so at a point with splitting number $0$.
This observation highlights the importance of the splitting number. 
It is not relevant to the case $n=1$ of course, but is a cornerstone of the
argument in higher dimensions.  

\begin{lemma}\label{le:split} Let $A_t$ be any path in $\G$ and
$\la(t)\in {\rm Spec}\,A_t$ a continuous family of eigenvalues  such that
\begin{eqnarray*} \la(t)\in S^1,\; \mbox{ for }\; t\le T,\\
\la(t)\notin S^1,\; \mbox{ for }\; t> T.
\end{eqnarray*}
Suppose also that $\la(t)$ has multiplicity $1$ when $t>T$ and 
multiplicity $2$ at $t = T$.  Then  
$$
\si(\la_{T}) = 0.
$$
\end{lemma}
\proof{}  For  $t\ge T$ let $V_t\subset \C^{2n}$ be the
space generated by the eigenspaces $E_{\la(t)}, E_{1/(\bar{\la}(t)}$.
By hypothesis, this is $2$-dimensional for each such $t$, and it
clearly varies continuously with $t$.  As above, $\be$ is
non-degenerate on each $V_t$.  Therefore the splitting number can be
$2, 0 $ or $-2$ and will be $0$ if and only if there are non-zero $v$
such that $\be(v,v) = 0$.  But when $t > T$, $\be(v,v) = 0$ on the
eigenvectors in $V_t$.  Therefore, $\be$ has signature $0$ when $t >
T$ and so must also have signature $0$ on $V_T$.
\QED

\section{Conjugacy classes and the positive cone}

\subsection{The case $n = 1$}

We will write $\Oo_{\Rr}^\pm$
for the set of all elements in $\Sl(2,\R)$  with 
eigenvalues in $\R^{\pm} - \{\pm 1\}$, and will also divide 
$\Uu - \{\pm \1, \pm\1^{\pm}\}$
into two sets $\Oo_\Uu^{\pm}$, distinguished by the splitting number of the
eigenvalue  $\la$ in the upper half-plane
$$
\Hh = \{ z\in \C: \Im\, z > 0\}.
$$
Thus, for elements of $\Oo_\Uu^+$, the eigenvalue in $\Hh$ has splitting number
$+1$, and, by considering the rotations
$$
\rho(t) = \left(\begin{array}{cc}  \cos t&-\sin t
\\ \sin t& \cos t \end{array}\right), 
$$
 it is easy to check that
$$
\Oo_\Uu^+ = \{A\,:\,  \tr A < 2,\; A_{21} > 0\},
\quad 
\Oo_\Uu^- = \{A\,:\,  \tr A < 2,\; A_{21} < 0\}.
$$
These sets project to the open strata in $\Co$.

We now work out in detail the structure of conjugacy classes near the 
element $-\1$.
The structure at $\1$ is similar, but will not be so important to us since 
we are interested
in paths which avoid the eigenvalue $1$.   

It is not hard to check that there are $3$ conjugacy classes with
spectrum $ \{-1\}$, the diagonalizable class consisting only of $\{-\1\}$, 
and the classes
$\Nn^\pm$ containing the nilpotent elements 
$$ 
\left(\begin{array}{cc}  -1&0
\\ \mp 1&-1 \end{array}\right). 
$$ 
Again, the classes $\Nn^{\pm}$ may be distinguished by the sign of the
$21$-entry.   Moreover, if we write the matrices near
$-\1$ as 
$$
 A = \left(\begin{array}{cc}  -1 + x & y \\ z& \frac{1 + yz}{-1+x}
\end{array}\right), 
$$
the elements in $\Nn^- \cup \{-\1\}\cup \Nn^+$ are those
with trace $= -2$ and so form the boundary of the cone  $x^2 + yz = 0$.  The
interior of the cone is the set where $x^2 < -yz $, and has two components
$$
\Uu^+ = \{y < 0, \,z > 0\},\quad \Uu^-=\{y > 0,\,z < 0\}.
$$ The rest of the space consists of points in $\Oo_{\Rr}^-$ with $x^2
> -yz $.  Observe that the labelling is chosen so that $\Nn^-$ is in
the closure of $\Oo_\Uu^+$.  (This class $\Nn^-$ was labelled by $-$
rather than $+$ because, as we shall see, it is the place where
positive paths leave $\Uu$.)

The next lemma describes the structure of short positive paths in $\SL(2,\R)$.

\begin{lemma}\label{le:pos1} Let $A_{t\in [0,1]}$ be a short positive path
in $\SL(2,\R)$, 
and for each $t>0$,  let
$$
\la_t = r(t) e^{i\theta_t},\quad \theta_t \in (0,2\pi),\;\;  |r(t)|\ge 1,
$$
be an eigenvalue  of 
$A_t$, chosen so that it has splitting number $+1$, when $\theta \ne \pi$.
Then:
\begin{description}
\item[(i)]
 there is $\eps > 0$ such that  $A_t \in \Oo_\Uu^+$ for $t \in [0,\eps]$;

\item[(ii)]  the function
$\theta_t$ is increasing, and is strictly increasing except perhaps at 
the point $\theta_t = \pi$;

\item[(iii)]   if $A_T = -\1$ for any $T$ then $A_t$ remains in $\Uu$ for 
$t > T$;

\item[(iv)] if $A_t$ enters $\Oo_{\Rr}^-$ it does so through a point of
$\Nn^-$ and if it reenters $\Uu$ it must go through  $\Nn^+$ to $\Oo_\Uu^-$.

\item[(v)]  $A_t\notin \Oo_{\Rr}^+$ for $t > 0$.
\end{description}
\end{lemma}

\proof{}  We first claim that at every point $N$ of $\Nn^-$ all vectors in the
positive tangent cone at $N$ point into $\Oo_\Rr^-$.  For if not, by
openness (Lemma~\ref{le:open}(i)), there would be a positive vector
pointing into $\Oo_\Uu^+$ and hence a positive path along which a
positive eigenvalue would move clockwise, in contradiction to Krein's
lemma.  Since conjugation by an orthogonal reflection interchanges
$\Nn^+$ and $\Nn^-$ and positive and negative, it follows that the
positive tangent cone points in $\Oo_\Uu^-$ at all points of $\Nn^+$.
Since every neighbourhood of $-\1$ contains points in $\Nn^\pm$, it
follows from openness that positive paths starting at $-\1$ must
remain in $\Uu$.  Since a similar result holds for negative paths,
there is no short positive path which starts with eigenvalues on
$\R^-$ and ends at $-\1$.  This proves (iii) and (iv).  A similar
argument proves (i) and (v).  Statement (ii) is a direct consequence
of Lemma~\ref{le:bott}).
\QED

\begin{remark}\label{rem:pm}\rm (i)  The above proof shows that the
structure of these positive paths in $\SL(2,\R)$ is determined by Krein's lemma
and the topology of the conjugacy classes.  However the result may
also be proved by direct computation.  For example, if 
$$
N = \left(\begin{array}{cc}  -1 & 0\\ -1 & -1
\end{array}\right)\in \Nn^+
$$
 and $P$ is
symmetric and positive definite, then
$$
\frac d{dt}\left|_{t=0}\right. {\rm tr\,}(e^{tJP}N)  = {\rm tr\,}(JPN) > 0.
$$ Hence the trace lies in $(-2,2)$ for $t> 0$, which means that every
positive path starting at $N$ moves into $\Uu^-$.  By invariance under
conjugacy, this has be to true for every element in $\Nn^+$.\MS

\NI
(ii)  It is implicit in the above proof that 
under a positive
flow in $\SL(2,\R)$ eigenvalues on $\R$ can flow in either direction.  
To see this explicitly, let  $A$ be the diagonal matrix $\diag(\la,1/\la)$ 
(where $|\la|> 1$)
and consider the positive matrices  
$$
P_1 = \left(\begin{array}{cc}  2 & -1\\
-1&1 \end{array}\right), \quad
P_2 = \left(\begin{array}{cc}  2 & 1\\
1&1 \end{array}\right).
$$
Then the derivative of
the trace (and hence the flow direction of $\la$) has different signs 
along the flows 
$t\mapsto e^{JP_it} A, i = 1,2$.
\end{remark}

  \subsection{The case $n=2$}

This case contains all the complications of the general case.  
We will see that the
behavior of  eigenvalues  along positive paths as they enter or 
leave the circle at a value in $S^1-\{\pm 1\}$ is much the
same as the behavior that we observed in the case $n = 1$ at the 
triple points $\pm 1$ in
$\Co$.    When $n = 2$ it is also possible
for two pairs of real eigenvalues to come together and then leave  
the real axis.  However,
we shall see that the positivity of the path imposes no essential 
restiction here.  One
indication of this is that there is no relevant notion of splitting 
number when $\la$ is not on the unit circle.

We first  consider the structure of elements $A\in \Sp(4,\R)$.  
As remarked above, the eigenvalues of $A$ consist either of
pairs $\la,\bar\la \in S^1$ and $\la, 1/\la \in \R - \{0\}$, or of 
quadruplets $\la,\bar\la, 1/\la,1/\bar\la$ where $\la\notin S^1\cup
\R$. We will label these pairs and quadruplets by the element
$\la$ with $|\la|\ge 1$ and $\Im\,\la \ge 0$.
\MS

A generic element of $\G$ lies in one of the following open regions:

\begin{description}
\item[(i)]  $\Oo_\Cc$, consisting of all matrices whose spectrum is a
quadruplet in $\C - (\R \cup S^1)$;
\item[(ii)]  $\Oo_\Uu$, consisting of all matrices whose eigenvalues 
 lie in $S^1 - \{\pm 1\}$ and either all have
multiplicity $1$ or have multiplicity $2$ and non-zero splitting 
numbers\footnote
{
By Lemma~\ref{le:split} the latter elements must be diagonalizable.
};
\item[(iii)]  $\Oo_{\Rr}$, consisting of all matrices whose eigenvalues have
multiplicity $1$ and lie in $\R-\{0,\pm 1\}$ 
(divided into two components: $\Oo_{\Rr}^{\pm}$);
\item[(iv)]  $\Oo_{\Uu,\Rr}$,  consisting of all matrices with $4$ distinct
eigenvalues, one pair on $S^1-\{\pm 1\}$ and the other on $\R-\{0,\pm 1\}$.
\end{description}

Note that  the  region $\Oo_{\Cc}$ is connected, but the others are not.

\begin{lemma}
The codimension $1$ part of the boundaries of the 
above regions
are:
\begin{description}  
\item[(v)]  $\Bb_\Uu$, consisting of all non-diagonalizable 
matrices whose spectrum consists of a
pair of conjugate points in $S^1 - \{\pm 1\}$ each of multiplicity $2$ and
splitting number $ 0$.  

\item[(vi)]  $\Bb_\Rr$, consisting of all non-diagonalizable
matrices whose spectrum is a
pair of distinct points $\la, 1/\la \in \R - \{\pm 1,0\}$ each of 
multiplicity $2$.

\item[(vii)]  $\Bb_{\Uu,\Rr}$, consisting of all non-diagonalizable matrices 
with spectrum
$\{\pm1, \pm1,$ $\la, \bar \la\}$, with $\la \in S^1 - \{\pm 1\}$.

\item[(viii)] $\Bb_{\Rr,\Uu}$, consisting of all non-diagonalizable matrices 
with spectrum
$\{\pm1, \pm1,$ $ \la, 1/\la\}$, with $\la \in \R - \{\pm 1,0\}$. 

\end{description}
\end{lemma}
\proof{}  Use Lemma~\ref{le:split}. 
\QED

 The next step is to  describe the conjugacy
classes  of elements in the above parts of  $\G$.  It follows easily from
Lemma~\ref{le:split} that the conjugacy class of an element $A\in \Oo_\Uu
\cup \Oo_{\Rr}\cup \Oo_{\Uu, \Rr}$ is determined by its spectrum and the
splitting numbers of those eigenvalues lying on $S^1$.  The next lemma
deals with the other cases.

\begin{lemma}\label{le:conj} \begin{description} \item[(i)]  The
conjugacy class of $A\in \Oo_\Cc$ is entirely determined by its eigenvalue
$\la \in \Hh$ with $|\la|> 1$.

\item[(ii)]   There are two types of conjugacy class in $\Bb_\Uu$. 
More precisely,  for each eigenvalue $\la \in S^1\cap \Hh$ of multiplicity $2$
and splitting number $0$, there are two conjugacy classes of non-diagonalizable
matrices, namely the two nilpotent classes $ \Nn_\la^\pm$.

\item[(iii)]  For
each $\la \in (-\infty, -1)\cup(1,\infty)$, there is one type of conjugacy 
class in $\Bb_\Rr$.

\item[(iv)] A matrix in $\Bb_{\Uu,\Rr}$ or $\Bb_{\Rr,\Uu}$ is conjugate
to a matrix which respects the splitting $\R^4 = \R^2 \oplus \R^2$. Thus the
conjugacy classes of $\Bb_{\Uu,\Rr}$ are determined by $\la \in (0,\pi)$, 
its splitting
number, and the conjugacy class $\Nn^-$ or $\Nn^+$ of the case $n=1$.
Those of $\Bb_{\Rr,\Uu}$ are 
determined by $\la \in (-\infty,-1) \cup (1, \infty)$ 
and the conjugacy class $\Nn^-$ or $\Nn^+$ of the case $n=1$. 
  \end{description}
\end{lemma}

\proof{}   First consider an element $A\in \Oo_\Cc$. By Lemma~\ref{le:be}(iv),
the subspace $V$ spanned by the eigenvectors $v,w$ with eigenvalues
$\la, 1/\bar\la$ is $\be$-orthogonal to its complex conjugate $\bar
V$, because this is the span of the eigenvectors $\bar v, \bar w$ with
eigenvalues $\bar\la, 1/\la$. By the same Lemma, $\be(v,v) = \be(w,w)
=0$.  Hence, because it is non-degenerate, the form $\be$ is non-zero
on $(v,w)$.  Thus $\be$ has zero signature on the subspaces $V, \bar
V$. Conversely, suppose given a $2$-dimensional subspace $V$ of $\C^4$
which is $\be$-orthogonal to its complex conjugate $\bar V$. Let $B$
be any complex linear $\be$-preserving automorphism of $V$ and extend
$B$ to $\bar V$ by complex conjugation.  Then $B$ is an element of
$\Sp(4,\R)$ because it is real and preserves $\be$.  In particular, if
$\be$ has signature $0$ on $V$, the subspace $V$ is spanned by vectors
$v,w$ such that $\be(v,v) = \be(w,w) = 0$.  Then $\be(v,w)
\ne 0$ since $\be$ is non-degenerate, and if we define $B = B(V,\la)$
by setting: 
$$
Bv = \la v,\quad Bw = (1/\bar\la) w,
$$
$B$ preserves $\be$ because
$$
\be(Bv,Bw) =  \be(\la v,1/\bar\la w) = \be(v,w).
$$
The corresponding element of $\Sp(4,\R)$ has spectrum
$\la,\bar\la,1/\la,1/\bar\la$.  Moreover, every element of 
$\Sp(4,\R)$ with such spectrum  has this form.
In particular, every such element is determined by the choice of
$\la$ with $|\la| > 1, \Im \,\la > 0$ and of the (ordered) decomposition
$V = \C v\oplus \C w$ of $V$.   Hence there is exactly one conjugacy
class for each such $\la$.  This proves (i).

To prove (ii), first observe that under the given hypotheses the 
eigenspace $V =
E_\la$ has the same properties as in (i).  Namely,  $\be|_V$ is a
non-degenerate form of signature $0$, and $V$ is $\be$-orthogonal to $\bar V$.
 Moreover, the restriction of $A$ to $V$ can
be any $\be$-invariant complex linear map with eigenvalue $\la$ of
multiplicity $2$ and splitting number $0$.  Suppose that
$A|_V$  has an eigenvector $v$ with  $\be(v,v) \ne 0$.  Then
there is a $\be$-orthogonal vector $w$ such that $\be(w,w) =
-\be(v,v)$, and because
$$
\be(w,w) = \be(Aw,Aw) = \be(\la w + \mu v, \la w + \mu v) = \be(w,w)
+ |\mu|^2\be(v,v)
$$
we must have $\mu = 0$.  Thus $A$ is diagonalizable (over $\C$) in this case.  
If no such $v$ exists, we may choose a basis $v,w$ for $V$ such that
$\be(v,v)=\be(w,w) = 0$ and so that $Av = \la v, Aw = \la w + \mu v$, for
some $\mu \ne 0$.  If we also set $\be(v,w) = i$, it is easy to see that
$\be(Aw,Aw) = \be(w,w)$ only if $\la\bar\mu$ is real.  
Given $v$, our choices have
determined $w$ up to a transformation of the form $w\mapsto w' = w +
\ka v$, where $\ka \in \R$ since $\be(w',w') = 0$.  It is easy to check
that if we make this alteration in $w$ then $\mu$ does not change.  On
the other hand, if we rescale $v,w$ replacing them by $ v/\ka,
\bar\ka\,w$ then $\mu$ changes  by the positive factor $|\ka|^2$. 
Hence we may assume that $\mu = \pm \la$.  This shows that $A|_V$ is
conjugate to exactly one of the  matrices
$$
D_\la = \left(\begin{array}{cc}  \la & 0\\0 & \la
\end{array}\right), \quad
N_\la^- = \left(\begin{array}{cc} \la & 0\\ -\la & \la
\end{array}\right),\quad
N_\la^+ = \left(\begin{array}{cc}  \la & 0\\ \la & \la
\end{array}\right).
$$

 To prove (iii), observe that the eigenspace $W$ of $\la$
(where $|\la| > 1$) is Lagrangian, and by suitably conjugating $A$ by an
element of $\G$ we may suppose that the other eigenspace is $JW$.   It then
follows that $A$ has the form $\la C \oplus \la^{-1} J^T(C^{-1})^TJ$ for some
$C\in \SL(2,\R)$ with trace $= 2$.    In fact, any  linear map of the form
$B\oplus J^T(B^{-1})^TJ$, where $B \in \GL(2,\R)$,  is  in $\G$.  It
follows easily that there are two conjugacy classes, one  in which $C$ is
diagonalizable, and one in which $C$ is not. 

\MS
\NI
(iv) $\,$ Finally, consider a matrix $A$ in $\Bb_{\Uu,\Rr}$ with a
pair of eigenvalues $\la, \bar{\la}$, with $\la \in S^1 \cap \Hh$
where $\Hh$ is the upper half plane, and the eigenvalue say $-1$, with
multiplicity $2$.  Assume say that the splitting number of $\la$ is
positive. Let $V$ be the eigenspace generated by the eigenvectors
$v,w$ corresponding to the eigenvalues $\la, \bar{\la}$, and $V'$ the
invariant subspace associated to the double eigenvalue $-1$. By
Lemma~\ref{le:be}, $v$ is $\be$\/-orthogonal to $V'$ and $w$, and
$\be(v,v) > 0$ by assumption. Similarly, $w$ is $\be$\/-orthogonal to
$V'$ and $v$, and $\be(w,w) < 0$.  Thus the restriction of $\be$ to
$V'$ is a also a non-degenerate symmetric Hermitian form with zero
signature. Because $A$ is non-diagonalizable, its restriction to $V'$
is non-diagonalizable.  Note that both eigenspaces $V$ and $V'$ are
invariant under conjugation, and therefore that the restrictions of
$A$ to each subspace is real.  This shows that we can identify $V
\oplus V'$ with $\C(\R^2) \oplus \C(\R^2)$ by a linear map $f$ which
respects the factors, preserves $\be$ and is real, and such that $f A
f^{-1}$ is the direct sum of a rotation on the first factor and of an
element in $\Nn^{\pm}$ on the other factor.

   The proof is similar in the other cases.   
 \QED

We next investigate the relationship between the positive cone field 
and the projection $\pi: \G \to \Co$.  Krein's Lemma~\ref{le:bott} 
shows that the movement of
eigenvalues along positive paths in $\Oo_\Uu$ and $\Oo_{\Uu,\Rr}$ is
constrained.   However, this is not so for
the other open strata.

\begin{lemma}\label{le:cone}  The projection $\pi: \G\to \Co$ maps 
the positive
cone $\Pp_A$ at $A\in \Oo_\Cc$ onto the tangent space to $\Co$ at $\pi(A)$. 
A similar statement is
true for $A\in \Oo_\Rr$.
\end{lemma}
\proof{}  Recall from Lemma~\ref{le:conj} that an
element $A$ of $\Oo_\Cc$  defines a unique
splitting $V\oplus \bar V$, where $\be|_V$ has signature $0$ and where $A|_V$
has eigenvalues $\la, 1/\bar\la$ for some $\la \in \{\Im\, z > 0\}$.  
Moreover,
we may choose the eigenvectors $v,w$ for $\la, 1/\bar\la$ respectively so that 
$$
\be(v,v) = \be(w,w) =0,\quad\be(v,w) = i.
$$
Therefore, if we fix such $v,w$ we get a unique representative of each
conjugacy class by varying $\la$.   

Most such splittings $V\oplus \bar V$ are not
$J$-invariant, and so it is hard to describe the positive flows at $A$.  
However, because the positive cone field on $\G$ is 
invariant under conjugation (see Lemma~\ref{le:conj0}), we only
need to prove this statement for one representative of 
each conjugacy class, and so we
choose one where this splitting is $J$-invariant.  For example,
if we set $v = 1/\sqrt2(1,0,i,0)$ and $w = -Jv$, all the
required identities are satisfied.  Moreover, the transformation
$X$ which takes the standard basis to $v,w,\bar v,\bar w$ is unitary.  
Therefore, any linear automorphism $B$ of $V$ which is represented
by a positive definite matrix with respect to the basis $v,w$ extends by
complex conjugation to a linear automorphism $B\oplus \bar B$ of
$\C^4$ which is also positive definite with respect to the usual basis. 
This means that there are positive tangent vectors to $\G$
at $A$ which preserve $V$ and restrict to $JB$ there.  
We are now essentially reduced to
 the $2$-dimensional case.  Just as in Lemma~\ref{le:pos1},  
if $B$ is a real positive
definite matrix with negative $21$-entry then the eigenvalue 
$\la$ moves along a ray
towards $S^1$ and if its $21$-entry is positive $\la$ moves 
along this ray away from
$S^1$.    Since $\pi$ maps the positive cone onto an open convex cone, 
this means that
the image has to be the whole tangent space to $\Co$.    

The argument for $A\in \Oo_\Rr$ is similar but easier.   \QED

The next step is to study generic positive paths through the various
parts of the codimens\-ion $1$ boundary components. We will see that
Krein's lemma forces the behavior of positive paths at the boundary
$\Bb_\Uu$ to mimic that that near $-\1$ in the case $n=1$, while the
boundary $\Bb_{\Rr}$ is essentially not seen by positive paths.  Our
results may be summarized in the following proposition.

\begin{prop}\label{prop:pos2}
\begin{description}
\item[(i)]  If $A \in \Nn_\la^-\subset \Bb_\Uu$, then $\pi$ maps
 the positive cone
$\Pp_A$ into the set of vectors at $\pi(A)$ which point into
$\pi(\Oo_\Cc)$.  Similarly, if $A \in \Nn_\la^+$, then $\pi_*(\Pp_A)$
consists of vectors pointing into $\pi(\Uu)$.
\item[(ii)]  If $A\in \Bb_\Rr$ then $\pi_*(\Pp_A)$ contains vectors which point
into $\pi(\Oo_\Cc)$ as well as vectors pointing into $\pi(\Oo_\Rr)$.
\item[(iii)]  If $A\in \Bb_{\Uu,\Rr} \cup \Bb_{\Rr,\Uu}$ then 
$A$ is conjugate to an element which preserves the splitting $\R^4 =
\R^2\oplus \R^2$, and  $\pi_*$ takes the subcone in $\Pp_A$ formed by
vectors which preserve this splitting onto the whole image $\pi_*(\Pp_A)$.  In
other words, positive paths through $\Bb_{\Uu,\Rr}$ and
$\Bb_{\Rr,\Uu}$ behave just like
paths in the product $\Sl(2,\R)\times \Sl(2,\R)$.
\end{description}
\end{prop}
\proof{} (i)
Suppose first that $A \in \Nn_\la^-$.  We claim that there is a 
neighbourhood of
$A\in \G$ whose intersection with $\Uu$ consists of elements which
have an eigenvalue $\la'$ with splitting number $+1$ which is close to $\la$
and to the right of it, i.e with ${\rm arg }\,\la' < {\rm arg}\,\la$.  
To see this, 
observe first that by conjugacy invariance it suffices to prove
this for a neighbourhood of $A$ in the subgroup
$$
\G_V = \{A'\in \G\,:\, A'(V) = V\}.
$$
If $A'\in \G_V$, then  $A'|_V$ may be
written as  $ \la' Y$ where $\la'$ is close to $\la$, and $Y$ has 
determinant $1$
and preserves $\be$.  If  $$ 
 Y=\left(\begin{array}{cc}  a & b\\c & d
\end{array}\right)
$$
with respect to the basis $v,w$, then we must have $a\bar c \in
\R$, $b\bar d\in \R$ and $a\bar d - c\bar b = 1$ in order to preserve
$\be$  and $\tr Y = 2$. This implies that $Y$ must be in
 $\SL(2,\R)$.  It is also close to the element $\frac 1{\la} A|_V \in
\Nn^- = \Nn_{-1}^-$. Therefore, the result follows  
from the corresponding result in the case
$n=1$: see the discussion just before Lemma~\ref{le:pos1}. 

Since  eigenvalues with positive splitting number flow anticlockwise
under a positive path (by Krein), every positive path through $A\in \Nn_\la^-$
must have all eigenvalues off $S^1$ for $t > T$ and so must enter $\Oo_\Cc$. 
Similarly, positive paths through any point of $\Nn_\la^+$ flow into $\Uu$.  As
in the case when $\la = -1$, this implies that positive paths through points of
$\Dd_\la$ remain in $\Uu$. 
\MS

\NI
(ii)  
Suppose first that $A$ in the closure $ \bar{\Bb}_\Rr$ of $\Bb_\Rr$
is diagonalizable,
say $A = \diag(\la,1/\la,\la,1/\la)$.  Then, as in
Remark~\ref{rem:pm}, there is a positive flow starting at $A$ which 
keeps $A$ diagonalizable and in $\Bb_\Dd$, while decreasing $|\la|$.  
Since any sufficiently
$C^1$-small perturbation of this flow is still positive, we can clearly 
flow positively from
$A$ into both regions $\Oo_\Rr$ and $\Oo_\Cc$.   

Next observe that if we fix $\la \in \R$ and consider only those elements of
$\Bb_\Rr$ with fixed $\la$-eigenspace $W$ and $1/\la$-eigenspace $W'$,
then these form a cone $C(\la)$ whose vertex is the diagonalizable element
$D_\la$.  (The structure here is just like that near $-\1$  which was 
discussed in
\S3.1.)  Therefore, if $v$ is a positive vector at $D_\la$ which points 
into $\Oo_\Cc$, there is
a nearby vector $v'$ at a nearby (non-diagonalizable) point $A$ which 
also points into
$\Oo_\Cc$.  Moreover, we may assume that  $v'$ is positive by the openness 
of the positive
cone.    Therefore one can move from $A$ into $\Oo_\Cc$, and similarly  
into $\Oo_\Rr$. 
Since there is only one conjugacy class of such $A$ for each $\la$, this 
completes the
proof. 
\MS

\NI
(iii) As shown in Lemma~\ref{le:pos1}, the behavior of positive paths
at the two conjugacy classes $-\1^\pm$ is dictated by Krein's lemma
and the topology of $\Sl(2,\R)$.  A similar argument applies
here. \QED

\begin{remark}\rm (i) One can construct explicit paths from $A\in \Bb_\Rr$
into $\Oo_\Cc$ and $\Oo_\Rr$ as follows.  Note that the matrices $A
\in \Bb_{\Rr}$ all satisfy the relation 
$$
\si_2 = \frac {\si_1^2}{4} + 2,
$$
where $\si_j$ denotes the $j$th symmetric function of the eigenvalues. 
Similarly, one can easily check that the matrices in $\Oo_{\Cc}$ 
satisfy the condition  $\si_2 >\frac {\si_1^2}{4} + 2$ while those
in $\Oo_{\Rr}$ satisfy $\si_2 < \frac {\si_1^2}{4} + 2$.  
Therefore,  to see where
the positive path  $A_t = (1 +  tJP + \dots)A$ goes as $t$ 
increases from $0$ we
just have to compare the derivatives
$$
\si_2' =\frac {d}{dt}\left|_{t=0}\right.\si_2,\quad 
\frac {\si_1\si_1'}{2}  =\frac {\si_1}{2} \frac
{d}{dt}\left|_{t=0}\right.\si_1. $$
For example, if we take
$$
A = \left(\begin{array}{cccc}  \la & 0 & 0 & 0\\
0& \frac 1{\la} & 0 & \al\\
-\frac{\al}{\la^2} & 0 & \la &0\\
0 &0 & 0 & \frac 1 {\la}\end{array}\right),\quad
P = \left(\begin{array}{cccc}  1 & 0 & 0 & 0\\
0& 5 & -2 & 0\\
0 & -2 & 1& 0\\
0 & 0 & 0 & 1\end{array}\right),
$$
then
$$
(1 + tJPA) - \mu \1 = \left(\begin{array}{cccc} 
 \la - \mu - \frac {2\al t}{\la^2}
& -\frac{5t}{\la} & 2t\la & -5\al t\\ 
t\la & \frac 1{\la}-\mu  & 0 & \al\\
-\frac{\al}{\la^2} & 0 & \la -\mu &-\frac{t}{\la}\\
-\frac{\al t}{\la^2} &-\frac{2t}{\la} & t\la & \frac 1 {\la} - \mu
-2\al t\end{array}\right).
$$
Thus
$$
\frac { \si_1\si_1' }{2} = -2\al(\frac 2{\la} + \la +
\frac 1{\la^3})
$$
while
$$
\si_2' = -2\al(2\la + \frac 2 {\la^3}).
$$
Therefore
$$
\si_2' - \frac { \si_1\si_1' }{2} = -2\al(\la + \frac 1{\la^3} - 
\frac 2{\la}) =
-2\al ({\la}^{1/2} - {\la}^{-3/2})^2, 
$$
and the path goes into $\Oo_\Cc$ if $\al < 0$ and into $\Oo_\Rr$ if $\al > 0$.
\MS 

\NI
(ii) Observe also that we did {\it not} 
show that  $\pi$ takes
the positive cone  at a point of $\Bb_\Rr$ onto the full tangent 
space of $\pi(\Oo_\Cc) =
\Hh$.  However, because of Lemma~\ref{le:open}(iii), all that 
matters to us is that
movement between these zones is possible.
\end{remark}

\section{Proof of the main results}

\subsection{The case $n \le 2$} 

For the convenience of the reader we now restate Theorem~\ref{thm:pospath}.

\begin{theorem}  Let $n \le 2$.
\begin{description}
\item[(i)]  An 
element of $\G - \Ss_1$ is the endpoint of  a short positive 
path (from $\1$) if and only if it has an even
number of real eigenvalues $\la$ with $\la > 1$.  
\item[(ii)]  There is a positive path between any two elements $A,B\in\G$.
Moreover, any short
positive path from $\1$ may be extended to a short positive path  with
endpoint in $\Uu$.
\item[(iii)]  The space of short positive paths (from $\1$) with endpoint 
in $\Uu$ is path-connected. 
\end{description}
\end{theorem}

\proof{}  (i) There are two ways eigenvalues can reach the positive real 
axis $\R^+$. Either
a pair of eigenvalues moves from $S^1$ through $+1$ to a pair on $\R^+$, 
or  a quadruplet
moves from $\C - \R$ to $\R$.  The first scenario cannot happen along a 
short path:  for
the only time a short path has an eigenvalue $+1$ is at time $t = 0$ and 
then, by
Lemma~\ref{le:pos1}, the eigenvalues have to move into $S^1$.  
In the second case, one
gets an even number of eigenvalues on $(1,\infty)$.
\MS

\NI
{(ii)}  Using Proposition~\ref{prop:pos2} it is not hard to see that 
there is a positive
path in $\Co$ from $\pi(A)$ to $\1$ and from $\1$ to $\pi(B)$.  
The result now follows from
Lemma~\ref{le:conj0}.  The second statement  is proved similarly.
\MS 

\NI
{(iii)}  We have to show that any  two short positive paths $\{A_t\},\{B_t\}$ 
beginning at $\1$ and ending in $\Uu$ can be joined
by a homotopy of short positive paths beginning at $\1$
and ending in $\Uu$.  We may clearly assume
that both paths are generic.  (Note however that the homotopy of paths may
have to go through codimension $2$ strata.) 
Moreover, by Lemmas~\ref{le:open} and~\ref{le:bott} 
 we may assume that there is $\delta > 0$
such that $A_t = B_t\in \Uu$ for $t\in [0,\delta]$.  Therefore, it
will suffice to consider the case when the second path $\{B_t\}$ is (a 
reparametrization of) $A_t, 0\le t \le \delta$.  Thus we have to show how to
\lq\lq shrink'' the path $\{A_t\}$ down to its initial segment in $\Uu$,
keeping its endpoint in $\Uu$.

We do this by constructing a homotopy $\ga^\mu_t$ in $\Co$ between the
 projections of the two paths which at each time obeys the
 restrictions stated in Proposition~\ref{prop:pos2}. Indeed, such a
 homotopy means the existence of a homotopy $\tilde{\ga}^\mu_t$ of
 short positive paths in $\G$ from $\1$ and ending in $\Uu$ which
 satisfies: $\tilde{\ga}^0_t = A_t$ for all $t$, $\tilde{\ga}^\mu_0 =
 \1$ and $\tilde{\ga}^\mu_1 \in \Uu$ for all $\mu$, and
 $\pi(\tilde{\ga}^1_t) = \pi(B_t)$.  By Proposition~\ref{prop:lift},
 there is a homotopy of paths (with starting points equal to $\1$ and
 endpoints free in the conjugacy class of $B_1$) between
 $\{\tilde{\ga}^1_t\}$ and $\{B_t\}$. The composition of these two
 homotopies is the desired path from $\{A_t\}$ to $\{B_t\}$.

We construct the homotopy $\ga^\mu$ in $\Co$ 
by means of a smooth map $r:\mu \mapsto
r^\mu$
from the interval $[0,1]$ to the space of short  positive paths in
$\Co$ as follows.\footnote
{
As will become clear, we are only interested in $r^\mu(t)$ on the time 
interval $t\in[\mu, 1]$.
} 
Suppose that we have found $r^\mu$ which satisfy:

\NI
(i) $r^\mu (\mu) = \pi(A_\mu)$ and  $r^\mu(1) \in\pi(\Uu)$; 

\NI
(ii) the paths $r^\mu$ are constant 
(with respect to $t$) for $\mu$ near $0$.     

\NI
Then it is not hard to check that we may take $\ga^\mu$ to be:
$$ 
\ga^\mu(t) = \begin{array} {lll}\pi(A_t) & \mbox{if}& t\le
\mu\\
                           r^\mu(t)& \mbox{if}& t\ge \mu.
\end{array}
$$

We will construct the $r^\mu$ backwards, starting at $\mu = 1$.
Intuitively, $r^\mu$ should be the simplest path from $\pi(A_\mu)$ to
$\pi(\Uu)$ (i.e. it should cross the stratum $\Co_1$ as few times as
possible), and we should think of the homotopy $\ga^\mu$ as shrinking
away the kinks in $\{A_t\}$ as $\mu$ decreases.  To be more precise,
let us define the {\bf complexity} $c(\ga)$ of a generic path $\ga =
\{C_t\}$ to be the number of times that $\pi(C_t)$ crosses $\Co_1$.
Then we will construct the $r^\mu$ so that $c(\ga^\mu)$ does not
increase as $\mu$ decreases.  Moreover it decreases by $2$ every time
that $A_\mu$ moves (as $\mu$ decreases) from $\Oo_\Rr$ to $\Oo_\Cc$,
or from $\Oo_{\Uu, \Rr}$ to $\Oo_\Uu$ or from $\Oo_\Cc$ to $\Oo_\Uu$.

To this end we assume that, except 
near places where $\pi(A_\mu)\in \Co_1$, the paths $r^\mu$ have the following
form:
\MS

\NI
$\bullet$  if $A_\mu\in\Uu$ then $r^\mu$ is constant;

\NI
$\bullet$  if $A_\mu\in \Oo_\Cc$ is labelled with $\la = se^{i\theta}$ then
$r_\mu$ goes down the ray $s'e^{i\theta}$ until it meets  the circle in
$\pi(\Nn_{e^{i\theta}}^+)$ and enters $\pi(\Uu)$; 

\NI
$\bullet$ if $A_\mu\in\Oo_\Rr$ then $r^\mu$ moves the two eigenvalue
pairs together, pushes them into $\pi(\Oo_\Cc)$ and then follows the previous
route to $\pi(\Uu)$;

\NI
$\bullet$ if $A_\mu\in\Oo_{\Uu,\Rr}$ then $r^\mu$ fixes the eigenvalues in
$S^1$ and moves the real ones through $\pi(\Nn_{-1}^+)$ to $S^1$.
\MS

We now describe how to extend $r$ over each type of crossing.  While
doing this, there will be times $\mu$ at which we will want to splice
different positive paths $r^\mu$ together.  More precisely, at these
$\mu_j$ we will have two different choices for the path $r^{\mu_j}$
from $A_{\mu_j}$ to $\Uu$.  But these choices will be homotopic (by a
homotopy which respects the endpoint conditions) and so, if we
reparametrize the path $A_t$ so that it stops at $t = \mu_j$ for a
little while, we may homotop from one choice to the other.  To be more
precise, we choose a nondecreasing function $\rho:[0,1]\to [0,1]$
which is bijective over all points except the $\mu_j$ and is such that
$\rho^{-1}(\mu_j)$ is an interval, and then change the relation
between $r^\mu,
\ga^\mu$ and $A_t$ by requiring that 
\MS

\NI
$\bullet$ $r^{\mu}(t)$ is defined for $t\in [\rho(\mu), 1]$ and
$
r^\mu(\rho(\mu)) = A_{\rho(\mu)}$,

\NI
$\bullet$ $ \ga^\mu(t) = A_{t}, t\le \rho(\mu).
$
 \MS

First let us consider how to handle crossings of the stratum
$\Bb_\Uu$.  Such crossings take place either at $\Nn_\la^+$ or
$\Nn_\la^-$.  However, because positive paths starting at $\Nn_\la^+$
point into $\Uu$ there is no problem extending $r$ smoothly over this
type of crossing.  Observe that at this crossing the complexity
remains unchanged.  The problem comes when $A_\mu\in
\Nn_\la^-$.  To deal 
with this, first let $C_t$ be a positive path  which starts in
$\Uu$, crosses $\Nn_\la^-$ at time $\mu$, goes a little into $\Oo_\Cc$
and then crosses back into $\Uu$ through $\Nn_\la^+$.  We may choose
this path so that it is homotopic through positive paths with fixed
endpoints to a positive path lying entirely in $\Uu$.  (In fact, we
could start with a suitable path in $\Uu$ and then perturb it.)  By
Lemma~\ref{le:open}, we may suppose that $\pi(A_t) = \pi(C_t) $ for
$t$ in some interval $(\mu - \eps, \mu + \eps)$ and so (using the
splicing technique described above) we may choose $r$ so that, for
some $\mu'\in (\mu,\mu + \eps)$, we have $r^{\mu'} (t) = \pi(C_t)$ for
$t \ge
\mu'$.  We now use this formula to extend $r$ over the interval $(\mu - \eps,
\mu')$.   Then $r^{\mu-\eps} = \{C_t\}$ starts and ends in $\Uu$ and  goes
only a little way outside $\Uu$.  We now reparametrize $\{A_t\}$
stopping it at time $t = \mu - \eps$ so that there is time first to
homotop $r^{\mu-\eps}= \{C_t\}$ with fixed endpoints to a path in
$\Uu$ and then to shrink it, fixing its first endpoint $A_{\mu-\eps}$,
to a constant path.  This completes the crossing.  Note that the
complexity does decrease by $2$.

A little thought shows that the same technique may be used to deal
with crossings of the other boundaries $\Bb_{\Uu,\Rr}, \Bb_{\Rr,\Uu}$
and $\Bb_\Rr$.  For example, if as $\mu$ decreases $A_\mu$ passes from
$\Oo_\Cc$ to $\Oo_\Rr$, one easily extends $r$ but does not change the
complexity of $\ga^\mu$.  On the other hand, if $A_\mu$ passes from
$\Oo_\Rr$ to $\Oo_\Cc$ one needs to take more trouble in extending $r$
but in exchange one decreases the complexity.
\QED

\subsection{The case $n > 2$}\label{ss:ge3}

As previously explained, the eigenvalues of an element of $\Sp(2n,
\R)$ form complex quadruplets or pairs on $S^1$ or $\R$.  It is easy
to check that the only singularities (or bifurcations) encountered by
a generic $1$-dimensional path $\ga = A_{t\in [0,1]}$ are those in
which
\begin{description}
\item[(a)]  a pair of eigenvalues on $S^1$ or $\R$ become equal to $\pm 1$;
\item[(b)]  
two pairs of eigenvalues coincide on $S^1$ or $\R$ and then move into
$\C$ (or conversely).  \end{description} In particular, two
quadruplets or three pairs do not coincide generically.  Therefore,
for each $t$, the space $\R^{2n}$ decomposes into a sum of $2$- and
$4$-dimensional eigenspaces $E_1(t)\oplus\dots\oplus E_k(t)$.
Moreover, the interval $[0,1]$ may be divided into a finite number of
subintervals over which this decomposition varies smoothly with $t$.
(The type of the decomposition may be different in the different
pieces.)  We claim that within each such piece the eigenvalue flow is
just the same as it would be if the decomposition were fixed rather
than varying.  The reason for this is that all the restrictions on the
eigenvalue flow are forced by Krein's lemma, which is valid for
arbitrary variation of decomposition, together with topological data
concerning $\Conj$ (i.e. topological information on the way the types
of conjugacy classes fit together.) With these remarks it is not hard
to adapt all the above arguments to the general case.  In particular,
the proof of (iii) is not essentially more difficult when $n > 2$
since in this case we reduce the complexity of the path as we procede
with the homotopy.  It is important to note that there still are
essentially unique ways of choosing the paths $r$ from the different
components of the top stratum $\Cc_0$ into $\pi(\Uu)$.  (\lq\lq
Essentially unique" means that the set of choices is connected.)  For
example, if $A^\mu$ were on a stratum with $3$ eigenvalue pairs on
$\R^-$, then one might combine two to make a quadruplet in $\C$ which
then moves down a ray to $S^1$ leaving one pair to move through
$\Nn_{-1}^+$ to $S^1$, or one might move all $3$ pairs directly down
$\R$ to $S^1$.  But these paths are homotopic (through positive paths
with endpoints on $\pi(\Uu)$) and so it is immaterial which choice we
make.  \QED

\section{Positive paths in Hamiltonian systems and in Hofer's geometry}

      We present here a brief outline of the various ways in which
positive paths intervene in the stability theory of Hamiltonian systems,
and in Hofer's geometry. Actually, they are also a crucial
ingredient of  the theory of closed
geodesics as developed by Bott in a series of papers (\cite{B1,B2,B3}). 
But since the
application of positive paths to closed geodesics (and in particular to
the computation of the index of iterates of a given closed geodesic) can
also be presented within the framework of positive Hamiltonian systems
via the geodesic flow (where Ekeland's formula appears as a generalization
of Bott's iteration formula), this application is, at least theoretically,
reducible to the following one.

\subsection{Periodic Hamiltonian systems}

   Many  Hamiltonian systems are given as periodic perturbations
of autonomous flows. When both the autonomous Hamiltonian and the perturbation
are quad\-ratic maps $\Qq_t$, the fundamental solution of the
periodic system
\begin{eqnarray*}
 \dot x(t)  & = & -J\nabla_X \Qq_t(x(t)),  \quad \mbox{where}\quad \Qq_t =
\Qq_{t+1}, \\
x(0) &= &\xi
\end{eqnarray*}
is a path of matrices $A_t \in \G$ as in \S1.  (In other words, the
trajectory $x(t)$ which starts at $\xi$ is $A_t(\xi)$.)  By the
periodicity of the generating Hamiltonian, it is clear that $A_{k+t} =
A_t A_1^k$ for all integers $k$ and real numbers $t \in [0,1)$.

\begin{definition}\rm The above periodic system is {\bf stable}
if all solutions $x(t)$ remain bounded for all times (in other words,
it is stable if
there is a constant $C$ such that $\| A_t \| \le C$ for all $t>0$). 
A matrix $A$ is {\bf stable} if there is $C$ such that $\|A^k\| \le C$ for all
positive integers $k$.
\end{definition}

   Clearly, the above system is stable if and only if $A=A_1$ is stable.
But it is easily seen that a matrix $A$ is stable if and only if all
its eigenvalues lie on the unit circle. Thus the stability of the 
periodic Hamiltonian system is determined by the spectrum of its time-$1$ flow
$A$.

\begin{definition} \rm The Hamiltonian system is {\bf strongly stable} if
any $C^2$-small periodic (and quadratic) perturbation of it remains stable. 
A matrix
$A$ is {\bf strongly stable} if it has a neighbourhood which
consists only of stable matrices.
\end{definition}

   Since a $C^2$-small perturbation $\Qq'_t$ of the Hamiltonian
$\Qq_t$ leads to a flow $A'_{t \in [0,1]}$ which is $C^1$-close (as a
path) to the unperturbed flow, the time-$1$ map $A' = A'(1)$ is
$C^0$-close to the time-$1$ map $A=A(1)$. Conversely any matrix $A'$
close enough to $A$ in the symplectic group is the time-$1$ map of a
$C^2$-small quadratic perturbation of $\Qq_{t \in [0,1]}$. This means
that a periodic Hamiltonian system is strongly stable if and only if
the time-$1$ map $A$ is strongly stable. Hence such a system is
strongly stable exactly when its time-$1$ flow belongs to the {\em
interior} of the set of symplectic matrices with spectrum in
$S^1$. But we have seen in the previous sections that this interior
consists of all matrices with spectrum in $S^1 - \{\pm 1\}$ and
maximal splitting numbers.  (More precisely, at each multiple
eigenvalue the absolute value of the splitting number must equal the
multiplicity of the eigenvalue.)  This is the Stability Theorem due to
Krein and Gelfand-Lidskii.

   Now consider the case $n=2$ and assume that the stable periodic
Hamiltonian $\Qq_t$ is negative and that its flow $\{A_{t}\}$ is short
(that is the index of $\{A_{t}\}$ is zero).  Further, let ${\cal NU}$
be the union of all {\em open} strata for which none of the
eigenvalues lie on $S^1$, and let $e_{A}(s)$ be the number of times
that the path $\{A_{t}\}$ enters ${\cal NU}$ and comes back to
$\Oo_{\Uu}$ during the time interval $(0, s]$.

\begin{lemma} For every short positive path  $A_{t \in [0,1]}$, beginning 
at $\1$, $e_{A}(1) \le 1$.
\end{lemma}
\proof{}  Because ${\cal NU}$ and $\Uu$ are open, it suffices to prove this 
when the path $\{A_{t}\}$ is generic.  If $e_A(1) = 0$ there is
nothing to prove.  So suppose that $e_A(s) = 1$ for some $s \le 1$.
Then Proposition~\ref{prop:pos2} implies that on the interval $[0, s]$
either $A_t$ moves from $\Oo_\Cc$ to $\Oo_\Uu$ passing through some
conjugacy class $\Nn_\la^+$ at some time $t_0 < s$ or $A_t$ moves from
$\Oo_\Rr^-$ into $\Oo_\Uu$ via matrices with all eigenvalues on
$(-\infty, -1)\cup S^1$.  In either case, it follows from Krein's
lemma~\ref{le:bott} that the first eigenvalue of $A_s$ that one
encounters when traversing $S^1$ anticlockwise from $\1$ has negative
splitting number.  Therefore, by Lemma~\ref{le:split}, the eigenvalues
of $A_t, t\ge s$, can only leave $S^1$ after a pair has crossed $+1$.
\QED

More generally, this argument shows that for a stable positive path
$A_{t \in [0,1]}$, which is not necessarily short or generic, $e_A(1)
\le i + 1$, where $i=i_A(1)$ is the Conley-Zehnder index of $A_{t \in
[0,1]}$.  (In this context, $i_A(1)$ is defined to be the number of
times $t > 0$ that a generic positive perturbation $\{A_t'\}$ of
$\{A_t\}$ (with $A_0' = A_0 = \1$) crosses the set $\Ss_1 = \{A: \det
(A - \1) = 0\}$.  Observe that because $\{A_t'\}$ is positive each
such crossing occurs with positive orientation.)

 Combining this with Ekeland's
version of Bott's iterated index formula (see \cite{EK}), we get:

\begin{proposition}\label{prop:restr} Let 
$A_{t \ge 0}$ be a stable positive path
generated by a $1$-periodic Hamiltonian in $\R^4$. Then there is a
positive real number $a$ such that for any time $t > 0$: $$ e_A(t) \le
a \, [t] \, i_A(1) + 1 $$ where $[t]$ is the greatest integer less
than or equal to $t$.
\end{proposition}

   Thus some aspects of the qualitative behavior of such paths are linearly
controlled by the Conley-Zehnder index. Similarly, the number of times that
such a path $A_t$ leaves $\Oo_{\Uu}$, enters any other open stratum
and comes back to $\Oo_{\Uu}$ during the time interval $(0,t]$,
is bounded above by
$$
3 (i_A(t) + 1) \le 3 (a \, [t] \, i_A(1) + 1).
$$

   As a last application, noted by Krein in \cite{GK}, Chap VI, consider
the system $\dot x(t) = -\mu J\nabla_X \Qq_t (x(t)), \, \Qq_{t+1} = \Qq_t$,
with a parameter $\mu \in \R$. Then:

\begin{prop}[Krein] If each $\Qq_t$ is negative definite,
there is $\mu_0 > 0$ such that the periodic Hamiltonian system 
$$
\dot x(t) = - \mu J \nabla_X \Qq_t (x(t)) 
$$
is strongly stable for all $0 < \mu < \mu_0$.
\end{prop}

    Actually, the value $\mu_0$ is the smallest $\mu$ such that the 
time-$1$ flow $A_{\mu}$ of the system has at least one eigenvalue 
equal to $-1$. Hence each time-$1$ flow $A_{\mu}$ for $\mu<\mu_0$
must be strongly stable.

\subsection{Hofer's geometry}

   Hofer's geometry is the geometry of the group $\Ham(M)$ of all 
smooth Hamiltonian
diffeomorphisms, generated by Hamiltonians $H_{t \in [0,1]}$ with
compact support, of a given symplectic manifold $M$, endowed with
the bi-invariant norm
$$
\|\phi\| = \inf_{H_t} \int_0^1 (\max_M H_t - \min_M H_t) \, dt
$$ where the infimum is taken over all smooth compactly supported
Hamiltonians $H_{t \in [0,1]}$ on $M$ whose time-$1$ flows are equal
to $\phi$. This norm defines a Finslerian metric (which is
$L^{\infty}$ with respect to space and $L^1$ with respect to time) on
the infinite dimensional Fr\'echet manifold $\Ham(M)$.  In \cite{LM1}
and \cite{BP}, it is shown that a path $\psi_t$ generated by a
Hamiltonian $H_t$ is a geodesic in this geometry exactly when there
are two points $p,P \in M$ such that $p$ is a minimum of $H_t$ and $P$
is a maximum of $H_t$ for all $t$. We then showed in \cite{LM2} that
the path $\psi_{t \in [0,1]}$ is a stable geodesic (that is a local
minimum of the length functional) if the linearized flows at $p$ and
$P$ have Conley-Zehnder index equal to $0$.  (See also~\cite{UST}.)
Note that the linear flow at $p$ is positive while the one at $P$ is
negative. The condition on the Conley-Zehnder index means that both
flows are short. The proof of the above sufficient condition for the
stability of geodesics in \cite{LM2} relies on the application of
holomorphic methods which reduce the problem to the purely topological
Main Theorem~\ref{thm:pospath}. Hence the result of the present paper
can be considered as the topological part of the proof of the
stability criterion in Hofer geometry. We refer the reader to
\cite{LM2} for the equivalence between this stability criterion and
the local squeezability of compact sets in cylinders.

\end{document}